\documentclass[12pt]{article}
\usepackage{amsfonts,amssymb,epsfig}
\usepackage{amsmath}

\date{}
\begin{document}
\newtheorem{df}{Definition}
\newtheorem{thm}{Theorem}
\newtheorem{lm}{Lemma}
\newtheorem{pr}{Proposition}
\newtheorem{co}{Corollary}
\newtheorem{re}{Remark}
\newtheorem{note}{Note}
\newtheorem{claim}{Claim}
\newtheorem{problem}{Problem}

\def\R{{\mathbb R}}

\def\E{\mathbb{E}}
\def\calF{{\cal F}}
\def\N{\mathbb{N}}
\def\calN{{\cal N}}
\def\calH{{\cal H}}
\def\n{\nu}
\def\a{\alpha}
\def\d{\delta}
\def\t{\theta}
\def\e{\varepsilon}
\def\t{\theta}
\def\pf{ \noindent {\bf Proof: \  }}
\def\trace{\rm trace}
\newcommand{\qed}{\hfill\vrule height6pt
width6pt depth0pt}
\def\endpf{\qed \medskip} \def\colon{{:}\;}
\setcounter{footnote}{0}

\def\Lip{{\rm Lip}}

\renewcommand{\qed}{\hfill\vrule height6pt  width6pt depth0pt}

\title{Matrix Subspaces of $L_1$
}

\author{Gideon Schechtman\thanks{Supported in part by the Israel Science Foundation. } } \maketitle

\begin{abstract}
If $E=\{e_i\}$ and $F=\{f_i\}$ are two 1-unconditional basic sequences in $L_1$ with $E$ $r$-concave and $F$ $p$-convex, for some $1\le r<p\le 2$, then the space of matrices $\{a_{i,j}\}$ with norm
$
\|\{a_{i,j}\}\|_{E(F)}=\big\|\sum_k \|\sum_l a_{k,l}f_l\|e_k\big\|
$
embeds into $L_1$. This generalizes a recent result of Prochno and Sch\"utt.
\end{abstract}

{\section{Introduction}

Recall that a basis $E=\{e_i\}_{i=1}^N$ of a finite ($N<\infty$) or infinite ($N=\infty$) dimensional real or complex Banach space is said to be $K$-unconditional if $\|\sum_i a_i e_i\|\le K\|\sum_i b_i e_i\|$ whenever $|a_i|=|b_i|$ for all $i$. Given a finite or infinite 1-unconditional basis, $E=\{e_i\}_{i=1}^N$, and a sequence of Banach spaces $\{X_i\}_{i=1}^N$ denote by $(\sum\bigoplus X_i)_E$ the space of sequences $x=(x_1,x_2,\dots)$, $x_i\in X_i$, for which the norm $\|x\|=\big\|\sum_i\|x_i\|e_i\big\|$ is finite.

If $X$ has a $1$-unconditional basis $F=\{f_j\}$ then $(\sum\bigoplus X)_E$ can be represented as a space of matrices $A=\{a_{i,j}\}$, denoted $E(F)$, with norm
\[
\|A\|_{E(F)}=\big\|\sum_i\|\sum_j a_{i,j}f_j\|e_i\big\|.
\]
In \cite{ps}, Prochno and Sch\"utt gave a sufficient condition for bases $E$ and $F$ of two Orlicz sequence spaces  which assure that $E(F)$ embeds into $L_1$. Here we generalize this result by giving a sufficient condition on two unconditional bases $E,F$, which assure that $E(F)$ embeds into $L_1$. As we shall see this condition is also ``almost" necessary.

Recall that an unconditional basis $\{e_i\}$ is said to be $p$-convex (resp. $r$-concave) with constant $K$ provided that for all $n$ and all $x_1,x_2,\dots,x_n$ in the span of $\{e_i\}$,
\[
\|\sum_{i=1}^n(|x_i|^{p})^{1/p}\|\le K (\sum_{i=1}^n\|x_i\|^{p})^{1/p}
\]
(resp.
\[
(\sum_{i=1}^n\|x_i\|^{r})^{1/r}\le K \|\sum_{i=1}^n(|x_i|^{r})^{1/r}\|\ ).
\]
Here, for $x=\sum x(j)e_j$ and a positive $\alpha$, $|x|^{\alpha}=\sum |x(j)|^{\alpha}e_j$.

$L_p$ will denote here $L_p([0,1],\lambda)$, $\lambda$ being the Lebesgue measure. As is known and quite easy to prove, any $1$-unconditional basic sequence in $L_p$, $1\le p\le 2$ (resp. $2\le p<\infty$), is $p$-convex (resp. $p$-concave) with constant depending only on $p$. It is also worthwhile to remind the reader that any $K$-unconditional basic sequence in $L_p$ is equivalent, with a constant depending only on $p$  and $K$ to a $1$-unconditional basic sequence in $L_p$.
It is due to Maurey \cite{ma} (see also \cite[III.H.10]{wo}), that for every $1\le r<p\le 2$, the span of every $p$-convex $1$-unconditional basic sequence in $L_1$ embeds into $L_p$ and also embeds into $L_r$ after change of density; i.e., there exists a probability measure $\mu$ on $[0,1]$ so that this span is isomorphic (with constants depending on $r,p$ and the $p$-convexity constant only) to a subspace of $L_r([0,1],\mu)$ on which the $L_r(\mu)$ and the $L_1(\mu)$ norms are equivalent.

If $M$ is an Orlicz function then the Orlicz space $\ell_M$ embeds into $L_p$ if and only if $M(t)/t^p$ is equivalent to an increasing function and $M(t)/t^2$ is equivalent to a decreasing function. This happens if and only if the natural basis of $\ell_M$ is $p$-convex and $2$-concave.

 Theorem \ref{thm:main} below states in particular that if $E$ and $F$ are two $1$-unconditional basic sequences in $L_1$ with $E$ $r$-concave and $F$ $p$-convex for some $1\le r<p\le 2$ then $E(F)$ embeds into $L_1$. When specializing to Orlicz spaces, this implies the main result of \cite{ps}.

\section{The main result}

\begin{thm}\label{thm:main}
Let $E=\{e_i\}$ be a 1-unconditional basic sequence in $L_1$ with $\{e_i\}$ $r$-concave with constant $K_1$ and let $X$ be a subspace of $L_{1}([0,1],\mu)$ for some probability measure $\mu$ satisfying $\|x\|_{L_{r}([0,1],\mu)}\le K_2 \|x\|_{L_{1}([0,1],\mu)}$ for some constant $K_2$ and all $x\in X$. Then $(\sum\bigoplus X)_{E}$ embeds into $L_1$ with a constant depending on $K_1,K_2$ and $r$ only.\\
Consequently, if $E=\{e_i\}$ and $F=\{f_i\}$ are two 1-unconditional basic sequences in $L_1$ with $E$ $r$-concave with constant $K_1$ and $F$ $p$-convex with constant $K_2$, for some $1\le r<p\le 2$, then the space of matrices $A=\{a_{k,l}\}$ with norm
\[
\|A\|_{E(F)}=\|\sum_k \|\sum_l a_{k,l}f_l\|e_k\|
\]
embeds into $L_1$ with a constant depending only on $r,p,K_1$ and $K_2$.
\end{thm}

\pf
The $p$-convexity of $\{f_i\}$ implies that after a change of density the $L_1$ and $L_r$ norms are equivalent on the span of $\{f_i\}$. See \cite{ma}. That is, there is a probability measure $\mu$ on $[0,1]$ and a constant $K_3$, depending only on $r,p$ and $K_2$ such that $\|\sum a_j\tilde f_j\|_{L_r([0,1],\mu)}\le  K_3\|\sum a_j\tilde f_j\|_{L_1([0,1],\mu)}$ for some sequence $\{\tilde f_j\}$ $1$-equivalent, in the relevant $L_1$ norm, to $\{f_j\}$, and for all coefficients $\{a_i\}$. It thus follows that the second part of the theorem follows from the first part.

To prove the first part, in $L_1([0,1]\times[0,1], \lambda\times\mu)$
consider the tensor product of the span of $\{e_i\}$ and $X$, that is the space of all functions of the form $\sum_i e_i\otimes x_i$, $x_i\in X$ for all $i$, where $e_i\otimes x_i(s,t)=e_i(s)x_i(t)$. Then, by the 1-unconditionality of $\{e_i\}$ and the triangle inequality,
\begin{eqnarray*}
\|\sum_i e_i\otimes x_i\|_1&=&\int\|\sum_i|x_i(t)|e_i\|_{L_1([0,1], \lambda)}d\mu(t)\\
&\ge&
\|\sum_i(\int|x_i(t)|d\mu(t))e_i\|_{L_1([0,1],\lambda)}\\
&=&\big\|\sum_i \|x_i\|e_i\big\|.
\end{eqnarray*}
On the other hand, by the 1-unconditionality and the $r$-concavity with constant $K_1$ of
$\{e_i\}$ (used in integral instead of summation form),
\begin{eqnarray*}
\|\sum_i e_i\otimes x_i\|_1&=&\int\int\big|\sum_i|x_i(t)|e_i(s)\big|d\lambda(s)d\mu(t)\\
&\le&
(\int(\int\big|\sum_i|x_i(t)|e_i(s)\big|d\lambda(s))^{r}d\mu(t))^{1/r}\\
&=&(\int\|\sum_i|x_i(t)|e_i\|_{L_1([0,1],\lambda)}^{r}d\mu(t))^{1/r}\\
&\le& K_1 \|\sum_i(\int|x_i(t)|^rd\mu(t))^{1/r}e_i\|_{L_1([0,1],\lambda)}\\
&\le& K_1K_2 \|\sum_i\int|x_i(t)|d\mu(t)e_i\|_{L_1([0,1],\lambda)}\\
&=&K_1K_2\big\|\sum_i \|x_i\|e_i\big\|
\end{eqnarray*}
\endpf

As is explained in the introduction the main result of \cite{ps} follows as corollary.
\begin{co}\label{co:ps}
If $M$ and $N$ are Orlicz functions such that $M(t)/t^r$ is equivalent to a decreasing function,
$N(t)/t^p$ is equivalent to an increasing function and $N(t)/t^2$ is equivalent to a 
decreasing function then $\ell_M(\ell_N)$ embeds into $L_1$.
\end{co}

\begin{re}
The role of $L_1$ in Theorem \ref{thm:main} can easily be replaced with $L_s$ for any $1\le s\le r$.
\end{re}
\begin{re}
If the bases $E$ and $F$ are infinite, say, and the smallest $r$ such that $E$ is $r$-concave is larger than the largest $p$ such that $F$ is $p$-convex, then $E(F)$ doesn't embed into $L_1$. This follows from the fact that in this case it is known that $\ell_r^n$ uniformly embed as blocks of $E$ and $\ell_p^n$ uniformly embed as blocks of $F$, for some $r>p$, while it is known that in this case  $\ell_r^n(\ell_p^n)$ do not uniformly embed into $L_1$.
\end{re}
This still leaves the case $r=p$, which is not covered in Theorem \ref{thm:main}, open:

\medskip\noindent
$\bullet$
{\em If $E$ and $F$ are two $1$-unconditional basic sequences in $L_1$ with $E$ $r$-concave and $F$ $r$-convex, does $E(F)$ embed into $L_1$?
}

\medskip

\noindent In the case that $E$ is an Orlicz space the problem above has a positive solution. We only sketch it. By the factorization theorem of Maurey mentioned above (\cite[III.H.10]{wo} is a good place to read it), and a simple compactness argument (to pass from the finite to the infinite case), it is enough to consider the case that $F$ is the $\ell_r$ unit vector basis. If the basis of $\ell_M$ is $r$-concave, then the $2/r$-convexification of $\ell_M$ (which is the space with norm $\|\{|a_i|^{2/r}\}\|_{\ell_M}^{r/2}$) embeds into $L_{2/r}$. This is again an Orlicz space, say, $\ell_{\tilde M}$. Now, tensoring with the Rademacher sequence (or a standard Gaussian sequence) we get that $\ell_{\tilde M}(\ell_2)$ embeds into $L_{2/r}$. We now want to $2/r$ concavify back, staying in $L_1$, so as to get that $\ell_M(\ell_r)$ embeds into $L_1$. This is known to be possible (and is buried somewhere in \cite{ms}): If $\{x_i\}$ is a $1$-unconditional basic sequence in $L_s$, $1<s\le2$ then its $s$-concavification (which is the space with norm $\|\{|a_i|^{1/s}\}\|_{\ell_M}^{s}$) embeds into $L_{1}$. Indeed, Let $\{f_i\}$ be a sequence of independent $2/s$ symmetric stable random variables normalized in $L_1$ and consider the span of the sequence \{$f_i\otimes |x_i|^s\}$ in $L_1$.

%
%

\noindent G. Schechtman\\
Department of Mathematics\\
Weizmann Institute of Science\\
Rehovot, Israel\\
{\tt gideon@weizmann.ac.il}
\\

\end{document}